# Efficient Solution Strategy for Chance-Constrained Optimal Power Flow based on FAST and Data-driven Convexification

Ren Hu
Dept. of Electrical & Computer Engineering
University of Central Florida
Orlando, USA

Qifeng Li
Dept. of Electrical & Computer Engineering
University of Central Florida
Orlando, USA

*Abstract*— The uncertainty of multiple power loads and renewable energy generations (PLREG) in power systems increases the complexity of power flow analysis for decision-makers. The chance-constrained method can be applied to model the optimization problems of power flow under uncertainty. This paper develops a novel solution approach for chance-constrained AC optimal power flow (CCACOPF) problem based on the data-driven convexification of power flow and a fast algorithm for scenario technique (FAST). This method is computationally effective for mainly two reasons. First, the original nonconvex AC power flow (ACPF) constraints are approximated by a set of learning-based quadratic convex ones. Second, FAST is an advanced scenario-based solution method (SSM) that doesn't rely on the pre-assumed probability distribution, using far less scenarios than the conventional SSM. Eventually, the CCACOPF is converted into a computationally tractable convex optimization problem. The simulation results on IEEE test cases indicate that 1) the proposed solution method can outperform the conventional SSM in computational efficiency, 2) the data-driven convexification of power flow is effective in approximating original complex AC power flow.

*Index Terms*-- Power flow, data-driven, chance-constrained, scenario optimization.

## I. INTRODUCTION

Solving optimal power flow is crucial to the decision-making on system planning, economic dispatching, reserve scheduling, etc., in power and energy systems. In practice, the renewable energy resources, such as wind and solar energy, are intermittent, and power load prediction carries unknown errors, both of which are the uncertainty sources considered in the optimization problems of power and energy systems, such as AC optimal power flow (ACOPF). Typically, there are stochastic [1], robust [2], and chance-constrained [3], [4] optimization approaches to model ACOPF with uncertainty. Unlike the stochastic method that pursues the expected objective value under a pre-assumed probability distribution and the robust method that concentrates on the worst case, the chance-constrained optimization (CCO) may interest the system decision-makers, providing a solution with a satisfactory confidence level (CL) to make sure that the system would run within the security boundary. The conventional solution methods for the CCO problems overly depend on the pre-assumed probability distribution (PD) of uncertain variables, which are called PD-based methods. However, such PD-based methods turn out to be computationally expensive for several reasons [5]-[7]. First, the multivariate integral of PD on the probabilistic constraints may be computationally prohibitive and cannot guarantee the convexity of constraints. Second, Monte Carlo simulation may be the only way to verify the solutions of the CCO problems, and it is too costly if a high CL is set. Instead of PD-based methods, the PD-free methods via random sampling (RS), called scenario optimization [5]-[7], only use the historical data of uncertainty and have attracted significant attention, widely used in probabilistic optimization problems. However, the number of samples required by the RS-based method (RSM) increases with the number of decision variables and may be very large [5], [6]. This makes RSM computationally ineffective. To tackle this computational issue, the fast algorithm for scenario technique (FAST) [7] is applied in this paper. It can greatly reduce the sample size for solving CCO problems and still enable the solution with a high CL. Besides, most applications of the scenario-based solution methods are under convex programs [5]-[7], but constraints of ACPF are originally nonconvex. The approximations of ACPF have been explored by either the linearization [4], [8], [9], or convexification techniques [10]-[16]. The linear methods, such as DC model [4] and linear ACPF [8], [9], may be easy to compute, but they ignore the quadratic terms of voltages, undermining the accuracy of model. Convex approximation methods such as the second-order cone (SOC) [10], semi-definite programming (SDP) [11], convex DistFlow [12], quadratic convex (QC) [13], moment-based [14], and the learning-based convex approximation [15], etc., have been frequently discussed. References [15], [16] have pointed out that the SDP-based relaxations cannot guarantee the exactness of solutions in some cases. The data-driven convex

This work is supported by the U.S. National Science Foundation under Grant #1808988.

quadratic approximation (DDCQA) [15] via ensemble learning exhibits better computational accuracy and efficiency than the SDP-based relaxation. Hence, the DDCQA of ACPF is adopted here. The contributions of this paper are summarized below: 1) A new solution strategy for CCACOPF problems is proposed using FAST with a small sample size and satisfactory CL; and 2) A data-driven convex approximation of power flow via ensemble learning is used to convert the intractable non-convex CCACOPF problem into a tractable convex optimization problem.

The rest of this paper is organized as follows. Section II illustrates the deterministic and chance-constrained ACOPF formulations. Section III discusses the solution approach for the CCACOPF problem via FAST and the DDCQA. Sections IV and V display the IEEE case analyses and conclusions, respectively.

## II. PROBLEM FORMULATION AND STATEMENT

This section: 1) formulates the deterministic ACOPF, 2) models the ACOPF under uncertainty using CCO approach is illustrated.

### A. Deterministic AC Optimal Power Flow

In an $n$-bus power system, the deterministic ACOPF can be presented by a series of equations and inequations in (1).

$$Min \sum_{i \in \vartheta}(a_{i0} + a_{i1}P_i^g + a_{i2}P_i^{g2}) \tag{1a}$$

$s.t.$

$$e_i \sum_{j=1}^n (G_{ij}e_j - B_{ij}f_j) + f_i \sum_{j=1}^n (G_{ij}f_j + B_{ij}e_j) = P_i^g - P_i^{net} \tag{1b}$$

$$f_i \sum_{j=1}^n (G_{ij}e_j - B_{ij}f_j) - e_i \sum_{j=1}^n (G_{ij}f_j + B_{ij}e_j) = Q_i^g - Q_i^{net} \tag{1c}$$

$$G_{ij}e_ie_j - B_{ij}e_if_j + B_{ij}f_ie_j + G_{ij}f_if_j - G_{ij}(e_i^2 + f_i^2) = P_{ij} \tag{1d}$$

$$G_{ij}f_ie_j - B_{ij}e_ie_j - G_{ij}e_if_j - B_{ij}f_if_j + B_{ij}(e_i^2 + f_i^2) = Q_{ij} \tag{1e}$$

$$V_i^{min2} \leq x_{2i-1}^2 + x_{2i}^2 \leq V_i^{max2} \tag{1f}$$

$$P_i^{gmin} \leq P_i^g \leq P_i^{gmax} \tag{1g}$$

$$Q_i^{gmin} \leq Q_i^g \leq Q_i^{gmax} \tag{1h}$$

$$P_{ij}^2 + Q_{ij}^2 \leq S_{ij}^{max} \tag{1i}$$

where $\vartheta$ is the index set of generators, ($i, j$ =1, 2, ..., $n$); $a_{i0}$, $a_{i1}$, $a_{i2}$ are the generator cost coefficients; $P_i^g$, $Q_i^g$ are the generator active and reactive power; $P_i^{net}$, $Q_i^{net}$ are the net active and reactive power inputs of power loads and renewable energy generations (PLREG); $P_i^{gmin}$, $P_i^{gmax}$, $Q_i^{gmin}$, $Q_i^{gmax}$ are the lower and upper limits of the generator active and reactive power; ; $V_i^{min}$, $V_i^{max}$ are the lower and upper limits of the bus voltage; $S_{ij}^{max}$ is the power line transmission capacity; $e_i$, $f_i$ represent the real and imaginary parts of voltage; $P_{ij}$, $Q_{ij}$, $G_{ij}$, $B_{ij}$ denote the active and reactive power of power lines, the real and imaginary parts of the line admittance.

### B. Chance-Constrained Formulation for ACOPF

Considering the uncertain net inputs of PLREG, the affine control scheme [17] is adopted to denote the fluctuating power generation as the base and adjustable parts, used to satisfy the forecast net inputs and the gap between the forecast and actual net inputs, respectively, shown in (2a)-(2d). We assume the nuance between the base part and real-time power loss is negligible and compensated by the reference generator.

$$P_{j,s}^{net} = P_{j(frt)}^{net} + \Delta p_{j,s} \tag{2a}$$

$$Q_{j,s}^{net} = Q_{j,t(frt)}^{net} + \Delta q_{j,s} \tag{2b}$$

$$P_i^g = P_{i(base)}^g - \omega_{pi} \sum_{j=1}^n \Delta p_{j,s}, \sum_{i \in \vartheta} \omega_{pi} = 1 \tag{2c}$$

$$Q_i^g = Q_{i(base)}^g - \omega_{qi} \sum_{j=1}^n \Delta q_{j,s}, \sum_{i \in \vartheta} \omega_{qi} = 1 \tag{2d}$$

where the index $s$ denotes under the $s$-th scenario; $\Delta p_{j,s}$, $\Delta q_{j,s}$ are the uncertain parts of $P_{j,s}^{net}$, $Q_{j,s}^{net}$; $P_{i(base)}^g$, $Q_{i(base)}^g$ are the base parts of $P_i^g$, $Q_i^g$; $\omega_{pi}$, $\omega_{qi}$ are the participation factors of generators. By substituting (2a)-(2d) into (1) and using the CCO, the deterministic problem (1) is converted into a CCO problem presented as:

$$\min l(y) \tag{2e}$$
$$s.t. \quad h(y, \delta) = 0 \tag{2f}$$
$$\Pr(g'(y, \delta) \leq 0) \geq 1 - \alpha \tag{2g}$$

where $h(y, \delta)$, $g'(y, \delta)$ compact constraints (1b)-(1e) and (1f)~(1i), respectively; $y$ is the decision variable vector, i.e., the active and reactive power of generators and power lines, and the bus voltage; $l(y)$ is the objective function; $\delta$ is the random variable vector, i.e., the net inputs of PLREG; Pr (·) denotes the probability that the constraint is satisfied; $\alpha$ is the CL factor.

## III. SCENARIO-BASED SOLUTION METHOD USING FAST AND DATA-DRIVEN CONVEXIFICATION OF POWER FLOW

This section introduces the RSM, i.e., scenario optimization, and then presents the proposed method for solving the CCACOPF problem via FAST and DDCQA. Eventually, the originally non-convex CCACOPF problem is converted to a convex optimization problem.

### A. Scenario Optimization

The conventional scenario-based method, i.e., RSM [5], [6], also called scenario optimization, plays a significant role in solving optimization problems with uncertainty, depending on the random samples. Under the convex optimization setting, the general formulation of scenario optimization is represented as

$$\min C^T z \tag{3a}$$
$$s.t. \quad F(z, \delta^s) \leq 0 \quad (s = 1,2,..,N') \tag{3b}$$

where (3a) is a linear objective function; $z$ is the decision variable vector; $F$ is convex in $z$; $\delta^s$ is the $s$-th random scenario of the random variable $\delta$; $N'$ is the required number of scenarios. More details of this method can be found in [5], [6]. To determine how many random samples are required to obtain the scenario-based solution of problem (3), the specific formulation can be written as

$$N' \geq \left(\frac{2}{\varepsilon}\right)\left(\ln\left(\frac{1}{\beta}\right) + d\right) \tag{3c}$$

where $d$ is the number of decision variables; $\varepsilon$ is the violation probability, $\varepsilon \in [0,1]$; $\beta$ is the CL factor. It is straightforward to see the minimum of sample size may explode as the number of decision variables increases in large-scale systems. This concern inspires us to further discuss the fast algorithm for scenario technique (FAST) in the next section.

### B. Fast Algorithm for Scenario Technique

Compared with the conventional scenario-based method, the fast algorithm for scenario technique (FAST) [7] is more efficient to solve the CCACOPF problem, due to its significant reduction of sample size and preserving the solution with a satisfactory CL. In FAST, the total number of samples is required at least N, consisting of $N_1$ and $N_2$, shown as

$$N_1 \geq d + 1 \tag{3d}$$
$$N_2 \geq \left(\frac{1}{\varepsilon}\right) \ln\left(\frac{1}{\beta}\right) \tag{3e}$$
$$N = N_1 + N_2 \geq \left(\frac{1}{\varepsilon}\right) \ln\left(\frac{1}{\beta}\right) + d + 1 \tag{3f}$$
$$\Pr(V(z^*) \leq \varepsilon) \geq 1 - \beta \tag{3g}$$

where $N_1$ is the initial number of scenarios; $N_2$ is the additional number of scenarios; $z^*$ is the solution provided by FAST. FAST [7] has proved that (3g) holds with $N$ samples smaller than the one obtained by RSM. The detailed procedure of FAST is summarized below in Algorithm 1.

---
**Algorithm 1: FAST**
1. Initialize the violation probability $\varepsilon$, the CL factor $\beta$, and the initial number of scenarios $N_1$ through (3d).
2. Compute the additional number of scenarios $N_2$ in (3h).
3. Solve the problem (3) with $N_1$ scenarios and set $(z_{N_1}^*, l_{N_1}^*)$ as the solution pair where $z_{N_1}^*$ is the output vector of generators; $l_{N_1}^*$ is the overall cost of generation in the objective function.
4. Detune the solution through $N_2$ scenarios by
$$z_{N_1+N_2}^* = \underset{1 \leq s \leq N_1+N_2}{\operatorname{argmax}} F(z_{N_1}^*, \delta^s) \tag{3h}$$
where $z_{N_1+N_2}^*$ is the output vector of generators updated by the additional $N_2$ new scenarios; $\delta^s$ is the $s$-th random scenario of the net inputs of PLREG. Note that (3h) is a simple linear program problem over $\delta^s$.
5. Output the objective cost $l_{N_1+N_2}^*$.

---

### C. Data-Driven Convex Quadratic Approximation of ACPF

Although FAST can mitigate the issue of large sample size in CCACOPF problem, constraints of ACPF are still non-convex. Therefore, the data-driven convex quadratic approximation (DDCQA) of power flow is applied [15]. For the ease of analysis, the constraints (1b)-(1e) are reformulated as the following matrix forms:

$$P_i = X^T A_i X \tag{4a}$$
$$Q_i = X^T B_i X \tag{4b}$$
$$P_{ij} = X_{ij}^T A_{ij} X_{ij} \tag{4c}$$
$$Q_{ij} = X_{ij}^T B_{ij} X_{ij} \tag{4d}$$
$$X = [x_1\, x_2, \ldots, x_{2n}]^T = [e_1\, f_1, \ldots, e_n\, f_n]^T \tag{4e}$$
$$X_{ij} = [x_{2i-1}\, x_{2i}\, x_{2j-1}\, x_{2j}]^T = [e_i\, f_i\, e_j\, f_j]^T \tag{4f}$$

where we uniformly redefine $x_i$ as the bus voltage variable. $A_i, B_i, A_{ij}, B_{ij}$ are symmetrical but indefinite matrices constituted by entries of the admittance matrix, indicating all dependent variables $P_i, Q_i, P_{ij}$, and $Q_{ij}$ in the power flow are non-convex functions of the independent variables $x_i$. Then, a convex quadratic mapping (4g)-(4j) between power, i.e., $P_i, Q_i, P_{ij}$, and $Q_{ij}$ and voltage $X$ or $X_{ij}$, is defined below:

$$P_i = X^T A_i^p X + B_i^p X + c_i^p \tag{4g}$$
$$Q_i = X^T A_i^q X + B_i^q X + c_i^q \tag{4h}$$
$$P_{ij} = X_{ij}^T A_{ij}^p X_{ij} + B_i^p X_{ij} + c_{ij}^p \tag{4i}$$
$$Q_{ij} = X_{ij}^T A_{ij}^q X_{ij} + B_i^q X_{ij} + c_{ij}^q \tag{4j}$$

where $A_i^p, A_i^q, A_{ij}^p, A_{ij}^q$ are positive semi-definite coefficient matrices of the quadratic terms, respectively; and $B_i^p, B_i^q, B_{ij}^p, B_{ij}^q$ are coefficient vectors of the linear terms; and $c_i^p, c_i^q, c_{ij}^p, c_{ij}^q$ are constant terms. To obtain the DDCQA of ACPF, polynomial regression (PR) is applied as a basic learner to learn the convex relationships between the bus voltage and the active or reactive power. Then, the ensemble learning method, i.e., gradient boosting is used to assemble all basic learners, to enhance the performance of model. The mean squared error function is used as the specific loss function

$$L(p_{mi}, \tilde{p}_{mi}) = \frac{1}{2}(p_{mi} - \tilde{p}_{mi})^2 \tag{4k}$$

where $p_{mi}$ and $\tilde{p}_{mi}$ are the observed and estimated values of $P_i$. The training dataset has $M$ samples $\{(X_m, Y_m)\}_{m=1}^M$ where $X_m$ is the $m$-th sample of bus voltages as the independent variables; $Y_m$ is the $m$-th sample of the active or reactive power. Taking $P_i$ as an example, the procedure of implementation is summarized in algorithm 2, named boosted PR.

---
**Algorithm 2: Boosted PR**
1. Initialize the model.
$$\boldsymbol{p}_{mi}^0(X) = \underset{\boldsymbol{\beta}'}{\arg\min} \sum_{m=1}^M L(p_{mi}, \boldsymbol{\beta}') \tag{4l}$$
where $\boldsymbol{\beta}'$ is the optimal constant.
2. For $\boldsymbol{t} = 1$: $T$ where $T$ is the maximum iterations.
  1) Compute the gradient descent $\boldsymbol{\gamma}_t$ by
$$\boldsymbol{\gamma}_t = -\left[\frac{\partial L(p_{mi}, p_{mi}(X))}{\partial p_{mi}(X)}\right]_{p_{mi}(X) = p_{mi}^{t-1}(X)} \tag{4m}$$
  2) Fit the PR model, $\boldsymbol{\varphi}_t(X; \boldsymbol{\theta})$ by
$$\boldsymbol{\theta}_t = \underset{\boldsymbol{\theta}}{\arg\min} \sum_{m=1}^M L(\boldsymbol{\gamma}_t, \boldsymbol{\varphi}_t(X_m; \boldsymbol{\theta})) \tag{4n}$$
where $\boldsymbol{\theta}_t$ denotes the regression coefficients, i.e., the learned parameters in equations (4a)-(4d).
  3) Update the model:
$$\boldsymbol{p}_m^t(X) = \boldsymbol{p}_m^{t-1}(X) + \boldsymbol{\mu}\boldsymbol{\varphi}_t(X_m; \boldsymbol{\theta}_t) \tag{4o}$$
where $\boldsymbol{\mu}$ is the learning rate.
3. Output $\boldsymbol{p}_m^T(X)$

---

### D. Scenario Programming Formulation for CCACOPF

Based on FAST and DDCQA, the new scenario-based solution approach for CCACOPF is written as a convex optimization problem below:

$$\min Z \tag{5a}$$
$$s.t. \quad (2e) \leq Z \tag{5b}$$
$$X_s^T A_i^p X_s + B_i^p X_s + c_i^p \leq P_i^g - P_{i,s}^{net} \tag{5c}$$
$$X_s^T A_i^q X_s + B_i^q X_s + c_i^q \leq Q_i^g - Q_{i,s}^{net} \tag{5d}$$
$$X_{ij,s}^T A_{ij}^p X_{ij,s} + B_i^p X_{ij,s} + c_{ij}^p \leq P_{ij,s} \tag{5e}$$
$$X_{ij,s}^T A_{ij}^q X_{ij,s} + B_i^q X_{ij,s} + c_{ij}^q \leq Q_{ij,s} \tag{5f}$$
$$x_{2i-1}^2 + x_{2i}^2 \leq V_i^{max2}, (s = 1,2,3\ldots,N) \tag{5g}$$
and constraints (1g)-(1i),

where $Z$ is the auxiliary variable to guarantee that the objective function is linear. The detailed implementation of solving problem (5) is provided in algorithm 1.

## IV. SIMULATION ANALYSIS

### A. Case Selection and Data Sampling

The real-world power systems with data are expected to be used in this research. However, they are not publicly available. As empirical alternatives, some IEEE standard test systems such as IEEE-5, -9, -57, -118 systems and the relevant simulation data are applied [15], [16]. The net active and reactive power loads at each bus are set up at the range of [0.7, 1.3] of their true values to simulate the uncertainty of PLREG For each test system the random scenarios are generated through Monte Carlo simulation. The goal is to demonstrate the efficacy of the proposed solution method for the CCACOPF problem with a small number of scenarios.

### B. Performance Comparison of Solution Approaches

For solving the CCACOPF problem, the violation probability and the CL factor are given as $\varepsilon=0.05$ and $\beta=10e-4$. The estimated number of scenarios required by FAST [7] and RSM [5], i.e., $N$ and $N'$ are shown in TABLE I below.

TABLE I
Comparison of Sample Size by FAST and RSM

| Size | Case5 | Case9 | Case57 | Case118 |
|---|---|---|---|---|
| $d$ | 32 | 42 | 240 | 716 |
| $N$ | 218 | 228 | 426 | 902 |
| $N'$ | 1649 | 2049 | 9969 | 29009 |
| $N'/N$ | 7.56 | 8.99 | 23.40 | 32.16 |

TABLE II
Comparison of Objective Costs by FAST and RSM (unit: $/hour)

| Cost | Case5 | Case9 | Case57 | Case118 |
|---|---|---|---|---|
| FAST_N1 | 24569.57 | 6872.37 | 14008.24 | 139194.50 |
| FAST_N2 | 26128.98 | 7345.04 | 14008.24 | 139194.50 |
| RSM | 26542.43 | 7396.21 | 14140.60 | 141988.90 |
| Difference | 1.56% | 0.69% | 0.94% | 1.97% |

TABLE I indicates that, compared to RSM, FAST greatly reduces the number of scenarios used for solving the CCACOPF problem. Especially in IEEE-118 system, the ratio of sample sizes $N'/N$ reaches up to more than 32 times. This implies that, in large-scale systems, FAST can work much faster than RSM. The optimal objective costs computed by FAST and RSM are displayed in TBALE II where $N_2 = 185$ is computed based on (3e). 'FAST_$N_1$' is the objective cost solved by FAST with $N_1$ scenarios. 'FAST_$N_2$' is the updated objective cost by FAST with $N_2$ scenarios. 'RSM' denotes the objective cost computed by RSM with $N'$ scenarios. 'Difference' is defined as (the difference between solutions of RSM and FAST) divided by the solution of RSM *100%. From TABLE II, we can observe that the objective cost updated with the additional $N_2$ new scenarios in FAST is almost same as the one computed by RSM. Specifically, the difference between solutions of RSM and FAST are bounded within 2% on all cases. In other words, FAST can not only retain the accuracy of the solution, but also use far less scenarios than RSM.

### C. Effect of the Violation Probability on the Solution

In practice, the system operators may prefer the solution with low violation probabilities. According to [5]-[7], the violation probability within 0.1 would give a satisfactory solution. To analyze the effect of the violation probability $\varepsilon$ on the solution, the sample sizes by FAST and RSM under $\varepsilon = 0.01\sim0.2$ and $\beta=10e-4$ are calculated and shown in TABLE III and IV, respectively. In FAST, as $N_1$ only depends on the number of decision variables for each case, changing $\varepsilon$ doesn't change $N_1$. And $N_2$ is inversely correlated to $\varepsilon$, shown in TABLE V.

TABLE III
Sample Sizes by FAST within $\varepsilon = 0.01\sim0.2$

| Size / $\varepsilon$ | 0.01 | 0.025 | 0.075 | 0.1 | 0.2 |
|---|---|---|---|---|---|
| Case5 | 955 | 402 | 156 | 126 | 80 |
| Case9 | 965 | 412 | 166 | 136 | 90 |
| Case57 | 1163 | 610 | 364 | 334 | 288 |
| Case118 | 1639 | 1086 | 840 | 810 | 764 |

TABLE IV
Sample Sizes by RSM within $\varepsilon = 0.01\sim0.2$

| Size / $\varepsilon$ | 0.01 | 0.025 | 0.075 | 0.1 | 0.2 |
|---|---|---|---|---|---|
| Case5 | 8243 | 3297 | 1099 | 825 | 413 |
| Case9 | 10243 | 4097 | 1366 | 1025 | 513 |
| Case57 | 49843 | 19937 | 6646 | 4985 | 2493 |
| Case118 | 145043 | 58017 | 19339 | 14505 | 7253 |

TABLE V
$N_2$ by FAST within $\varepsilon = 0.01\sim0.2$

| Size / $\varepsilon$ | 0.01 | 0.025 | 0.05 | 0.075 | 0.1 | 0.2 |
|---|---|---|---|---|---|---|
| $N_2$ | 922 | 369 | 185 | 123 | 93 | 47 |

TABLE VI
Objective Costs by FAST within $\varepsilon = 0.01\sim0.2$ (unit: $/ hour)

| Cost / $\varepsilon$ | 0.01 | 0.025 | 0.075 | 0.1 | 0.2 |
|---|---|---|---|---|---|
| Case5 | 26542.43 | 26293.85 | 25991.92 | 25991.92 | 25991.92 |
| Case9 | 7383.77 | 7383.77 | 7210.84 | 7210.84 | 7210.84 |
| Case57 | 14008.24 | 14008.24 | 14008.24 | 14008.24 | 14008.24 |
| Case118 | 139194.50 | 139194.50 | 139194.50 | 139194.50 | 139194.50 |

TABLE VII
Objective Costs by RSM within $\varepsilon = 0.01\sim0.2$ (unit: $/ hour)

| Cost / $\varepsilon$ | 0.01 | 0.025 | 0.075 | 0.1 | 0.2 |
|---|---|---|---|---|---|
| Case5 | 27047.25 | 26542.43 | 26542.43 | 26542.43 | 26293.85 |
| Case9 | 7567.79 | 7535.09 | 7383.77 | 7383.77 | 7383.77 |
| Case57 | 14263.63 | 13981.36 | 14140.60 | 14008.24 | 14008.24 |
| Case118 | 142679.58 | 141988.87 | 141988.87 | 141988.87 | 140799.73 |

TABLE III and IV show that: 1) the smaller $\varepsilon$, the larger the sample sizes by RSM and FAST, 2) the sample size by FAST is far smaller than the one by RSM. Meanwhile, TABLE V reveals that the increase in the sample size by FAST is only related to $N_2$ when shrinking $\varepsilon$, rather than the number of decision variables. In addition, the corresponding objective costs by FAST and RSM are computed in TABLE VI and VII, respectively. TABLE VI and VII indicate that with $\varepsilon=0.01$~$0.2$, the solutions provided by RSM or FAST have no significant changes. For example, on case 57 and case 118, the solutions provided by FAST don't change and the solutions provided by RSM have less than 2% differences. Similarly, on case 5 and case 9, the maximum difference between solutions provided by FAST or RSM is 2.34% or 3.40%. Eventually, the difference between solutions of RSM and FAST under different $\varepsilon$ is summarized in TABLE VIII and Fig. 1 below. In a word, under different $\varepsilon$, the difference between solutions of RSM and FAST is no more than 2.5%, i.e., FAST is a relatively robust alternative method for RSM to solve CCO problems.

TABLE VIII
Difference between Solutions of RSM and FAST within $\varepsilon = 0.01$~$0.2$

| Difference / $\varepsilon$ | 0.01 | 0.025 | 0.05 | 0.075 | 0.1 | 0.2 |
|---|---|---|---|---|---|---|
| Case5 | 1.87% | 0.94% | 1.56% | 2.07% | 2.07% | 1.15% |
| Case9 | 2.43% | 2.01% | 0.69% | 2.34% | 2.34% | 2.34% |
| Case57 | 1.79% | -0.19% | 0.94% | 0.94% | 0.00% | 0.00% |
| Case118 | 2.44% | 1.97% | 1.97% | 1.97% | 1.97% | 1.14% |

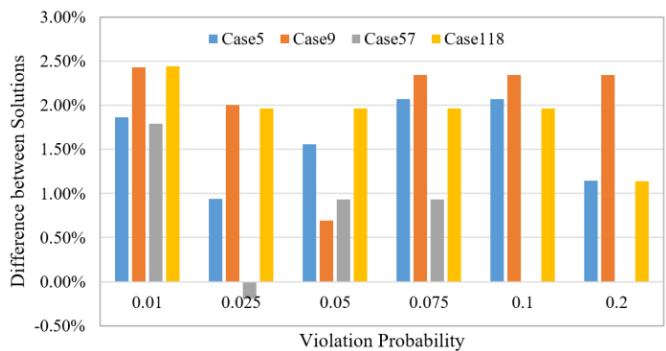

Fig. 1 Summary of Difference between Solutions

## V. CONCLUSIONS AND FUTURE WORK

This paper provides a novel approach to solve the intrinsically non-convex CCACOPF, based on the DDCQA of power flow and FAST. The proposed method only uses few scenarios without any assumed PD of the uncertain PLREG. As an improvement of the conventional scenario-based method, i.e., random sampling (RS)-based method (RSM), FAST cannot only greatly lessen the sample size required for solving the CCACOPF, but also maintain the accuracy of solution with a high CI. In other words, FAST is computationally faster than RSM, especially for large-scale systems with numerous decision variables. Basically, the scenario-based solution methods are applied mostly in the convex programs, but the original constraints of ACPF are non-convex. This gives rise to the necessity of approximating ACPF by convexification techniques. The DDCQA therefore is formulated using ensemble learning-based polynomial regression and applied to address the non-convex constraints of ACPF. Eventually, the computationally intractable CCACOPF is converted to a tractable convex quadratic optimization problem via the DDCQA and FAST. In our future work, we intend to use more real-life large-scale power systems to test the proposed solution method.